\documentclass[10pt]{amsart}
\UseRawInputEncoding
\usepackage{amsfonts}
\usepackage{amssymb}
\usepackage{amsmath,amsthm}
\usepackage{amscd}
\usepackage{graphics}
\usepackage{amsmath,amssymb}
\usepackage{pdfsync}
\usepackage{mathrsfs}
\usepackage{stmaryrd}
\usepackage[colorlinks,linkcolor=blue,citecolor=blue,urlcolor=blue]{hyperref}

\newcommand{\beq}{\begin{equation}}
\newcommand{\eeq}{\end{equation}}
\newcommand{\bea}{\begin{eqnarray}}
\newcommand{\eea}{\end{eqnarray}}

\newcommand{\bk}{\begin{cases}}
\newcommand{\ek}{\end{cases}}

\newtheorem{theorem}{Theorem}

\newtheorem{definition}[theorem]{Definition}

\newtheorem{lemma}{Lemma}

\newtheorem{remark}{Remark}

\numberwithin{equation}{section}
\begin{document}

\author{Piergiulio Tempesta}
\address{Departamento de F\'{\i}sica Te\'{o}rica, Facultad de F\'{\i}sicas, Universidad
Complutense, 28040 -- Madrid, Spain
 and Instituto de Ciencias Matem\'aticas, C/ Nicol\'as Cabrera, No 13--15, 28049 Madrid, Spain.}
\email{p.tempesta@fis.ucm.es, piergiulio.tempesta@icmat.es}
\title{On Appell sequences of polynomials of Bernoulli and Euler type}

\date{February 14, 2021}

\begin{abstract}
A construction of new sequences of generalized Bernoulli polynomials of
first and second kind is proposed. These sequences share with the classical Bernoulli polynomials many algebraic and
number--theoretical properties. A
new class of Euler--type polynomials is also presented.
\end{abstract}

\maketitle

\tableofcontents

\section{Introduction}

In this article \footnote{This is a revised version of the published paper. Some remarks have been added and misprints corrected.} several new classes of polynomials of Appell type, which
generalize the classical Bernoulli and Euler polynomials, will be introduced
and discussed.

We recall that the Bernoulli polynomials \cite{Bernoulli} are defined by the
generating function
\begin{equation}
\frac{t}{e^{t}-1}e^{xt}=\sum_{k=0}^{\infty }\frac{\,B_{k}\left( x\right) }{k!%
}t^{k}\text{.}  \label{I.1}
\end{equation}%
For $x=0$, formula (\ref{I.1}) reduces to the generating function of the
Bernoulli numbers:
\begin{equation}
\frac{t}{e^{t}-1}=\sum_{k=0}^{\infty }\frac{\,B_{k}}{k!}t^{k}\text{.}
\label{I.2}
\end{equation}%
The Bernoulli polynomials are also determined by the two properties
\begin{equation}
D\,B_{n}\left( x\right) =nB_{n-1}\left( x\right) \text{,}  \label{I.3}
\end{equation}%
\begin{equation}
\Delta \,B_{n}\left( x\right) =nx^{n-1}\text{,}  \label{I.4}
\end{equation}%
with the condition $B_{0}\left( x\right) =1$.\thinspace Here $D\,$is the
continuous derivative, $\Delta :=\left( T-1\right) $ is a discrete derivative,
and $T\,$ is the shift operator, defined by $Tf\left( x\right) =f\left(
x+1\right) $. Standard references are \cite{DSS} and \cite{Slav1} and the
books \thinspace \cite{Erdelyi} and \cite{IR}. By virtue of relation (\ref%
{I.3}), the Bernoulli polynomials belong to the class of Appell polynomials.

The Bernoulli polynomials play an important r\^{o}le in several branches of
mathematical analysis, such as the theory of distributions in $p$--adic
analysis \cite{Koblitz}, the theory of modular forms \cite{Lang2}, the study
of polynomial expansions of analytic functions \cite{BB}, etc. Formula (\ref%
{I.4}) is the basis for the application of these polynomials to interpolation theory \cite{Phil}. 
Several
generalizations of Bernoulli polynomials have also been proposed.
Particularly important is that of Leopold and Iwasawa, motivated by a
connection with the theory of $p$--adic L--functions \cite{Iwa}. Another
example is provided by the work of Carlitz \cite{Carl}. Recently, new
applications of the Bernoulli polynomials have also been found in
mathematical physics, in connection with the theory of the Korteweg--de
Vries equation \cite{FV} and Lam\'{e} equation \cite{GV}, and in the study
of vertex algebras \cite{DLM}.

The Bernoulli numbers are relevant in several branches of number theory, in
particular in the computation of special values of zeta functions (see, e.g.,  \cite{Eie}, \cite{Zag}), in the theory of cyclotomic fields \cite{Wash} and, since
Kummer's work, in connection with Fermat's last theorem \cite{IR}. They are
also useful in singularity theory \cite{Arnold1} and in the study of Coxeter groups \cite{Arnold2}. Standard applications to algebraic topology
are found in the computation of Todd characteristic classes and in the
Hirzebruch signature theorem, as well as, more recently, in complex homology
theory \cite{BCRS}, \cite{Ray}. In the last years, Bernoulli number
identities have found applications in quantum field theory \cite{DS} and in
the computation of Gromov--Witten invariants \cite{FP}.

In the papers \cite{Tempesta1}, \cite{Tempesta2} the following
generalization of the Bernoulli polynomials was introduced.

\begin{definition}
Let us consider the polynomial ring $\mathbb{Q}[c_{1},c_{2},...]$ and the formal power series
\begin{equation}
F\left( s\right) =s+c_{1}\frac{s^{2}}{2}+c_{2}\frac{s^{3}}{3}+ \ldots
\label{I.5}
\end{equation}
Let $G\left( t\right) $ be the compositional inverse series:
\begin{equation}
G\left( t\right) =t-c_{1}\frac{t^{2}}{2}+\left( 3c_{1}^{2}-2c_{2}\right)
\frac{t^{3}}{6}+\ldots  \label{I.6}
\end{equation}
so that $F\left( G\left( t\right) \right) =t$. The \textit{universal
higher--order Bernoulli polynomials }

$B_{k,a}^{G}\left(x,c_{1},\ldots,c_{n}\right) \equiv B_{k,a}^{G}\left(
x\right) $ are defined by
\begin{equation}
\left( \frac{t}{G\left( t\right) }\right) ^{a}e^{xt}=\sum_{k\geq
0}B_{k,a}^{G}\left( x\right) \frac{t^{k}}{k!}\text{,}  \label{I.7}
\end{equation}
where $a\neq 0$.
\end{definition}

This definition is clearly motivated by the works by Clarke \cite{Clarke},
Ray \cite{Ray} and Adelberg \cite{Adel} on \textit{universal Bernoulli
numbers}. We observe that many known polynomial sequences can be obtained by suitable choices of the coefficients appearing in the generating function \eqref{I.7}. For instance, the case of the standard Bernoulli polynomials corresponds to the choice $a=1$, $c_{i}=\left( -1\right)
^{i}$, since now $F\left( s\right) =\log \left( 1+s\right) \,$and $G\left(
t\right) =e^{t}-1$. When $a$ is a rational integer, $c_{i}=\left( -1\right)
^{i}\,$ we obtain the higher--order Bernoulli polynomials, which have also
been extensively studied (e.g. in  \cite{Carlitz}, \cite{CO},
\cite{How}, \cite{Roman}). When $x=0,$ $a\in \mathbb{R}$ , formula (\ref{I.7}) can be reduced to the generating function of the N\"{o}rlund polynomials and their generalizations.

By construction, the polynomials $B_{k,1}^{U}\left( 0\right) \equiv \widehat{%
B_{k}}\in \mathbb{Q}\left[ c_{1},\ldots,c_{n}\right] $, where $\widehat{B_{k}%
}\, $ are the universal Bernoulli numbers introduced by Clarke in \cite%
{Clarke}. The connection with algebraic topology relies on the fact that $%
G\left( F\left( s_{1}\right) +F\left( s_{2}\right) \right) \,$represents the
\textit{universal formal group} \cite{Haze}; the series $F$ is called the
formal group logarithm and $G$ the formal group exponential. The universal
formal group is defined over the Lazard ring $L$, which is the
subring of $\mathbb{Q}\left[ c_{1},\ldots,c_{n}\right] $ generated by the
coefficients of the power series $G\left( F\left( s_{1}\right) +F\left(
s_{2}\right) \right) $. In algebraic topology, in particular in complex
cobordism theory, the coefficients $c_{n}$ are identified with the cobordism
classes of $\mathbb{C}P^{n}$ (see \cite{AB}, \cite{Mil} and \cite{Ray}).
Generalizations of famous congruences valid for the classical Bernoulli
numbers, like the celebrated Kummer and Clausen--von Staudt congruences \cite%
{Adel}, \cite{IR}, are satisfied by these universal numbers.

The universal polynomials (\ref{I.7}) possess as well several interesting
general properties studied in \cite{Tempesta1}, \cite{Tempesta2}. 


The first aim of this paper is to provide nontrivial realizations of polynomial
sequences of the type (\ref{I.7}). In particular, we will focus on a
specific class of Bernoulli--type polynomials, constructed using the finite
operator calculus, as formulated by G. C. Rota and S. Roman \cite{Roman},
\cite{Rota}. For a recent review of the vast literature existing on this
approach (also known, in its earlier formulations, as the Umbral Calculus)
see, e.g., \cite{BL}. In other words, the \textit{Bernoulli--type
polynomials of order p} considered here correspond to a class of formal
power series $G\left( t\right) $ representing suitable \textit{difference delta
operators of order }$p$ (denoted by $\Delta _{p}$ ). These delta operators
have been introduced in \cite{LTW1}, \cite{LTW2}, where a version of Rota's
operator approach, based on the theory of representations of the
Heisenberg--Weyl algebra, has been outlined. In this paper we wish to
further illustrate the connection between Rota's approach and the theory of
Appell and Sheffer polynomials. The basic sequences associated with the
delta operators $\Delta _{p}$ are here studied in detail and some of their
combinatorial properties derived. In addition, their connection constants
with the basic sequence $\left\{ x^{n}\right\} $ (and the associated inverse
relations) define generalized Stirling numbers of
the first and second kind. The problem of classifying all \textit{Sheffer
sequences} associated with the delta operators considered in this paper is
essentially open.

By analogy with the classical case, the Bernoulli--type polynomials, parametrized by a real variable $a$ introduced in this work are  determined by the
relations
\begin{equation}
D\,B_{n,a}^{p}\left( x\right) =nB_{n-1,a}^{p}\left( x\right) \text{,}
\end{equation}%
\begin{equation}
\Delta^{a}_{p}\,B_{n,a}^{p}\left( x\right) =nx^{n-1}\text{,}
\end{equation}%
 and by the condition $B_{0,a}^{p}\left( x\right) =1$; here $\Delta^{a}_p$ is a suitable operator (defined in Sec. 5.2). For $a=1$, $p=1$, the usual characterization of the classical Bernoulli polynomials is recovered.
 

The use of the Roman--Rota formalism has the advantage of providing us with a
natural language which allows to unify the mathematical treatment of many polynomial
sequences.

The second aim of this paper is to introduce the notions of \textit{%
Bernoulli polynomials of the second kind} and of \textit{Euler--type
polynomials}.

All these polynomials can also play a role in discrete mathematics. Indeed,
they satisfy certain interesting linear difference equations of order $p$ in
one variable, defined in a two--dimensional space of parameters.

The future research plans on generalized Bernoulli polynomials include
possible applications in analytic number theory, in particular
to the theory of Dirichlet L--series and Riemann--Hurwitz zeta function.
Results in this direction have been obtained in \cite{Tempesta1} and \cite{Tempesta2}, where a
new construction relating generalized Bernoulli polynomials, formal groups
and L--series is proposed. 


The paper is organized as follows. In Section 2, a brief introduction to
finite operator calculus is presented, with a discussion in Section 3 of
difference delta operators. In Section 4, generalized Stirling numbers are
considered. In Section 5, new Bernoulli--type polynomials of the first kind
are introduced. The Bernoulli--type polynomials of the second kind are
introduced and discussed in Section 6. A generalization of Euler polynomials
is discussed in the final Section 7.

\section{Finite operator calculus}

In this section, some basic results concerning the theory of finite difference
operators and their relation with polynomial sequences of Appell and Sheffer tyoe are
reviewed. For further details and proofs, see the monographs \cite{Roman},
\cite{Rota}, where an extensive and modern treatment of this topic is
proposed.

Let $\mathcal{F}$ denote the algebra of formal power series in one variable $%
t$, endowed with the operations of sum and multiplication of series. Let $%
\mathcal{P}$ be the algebra of polynomials in one variable $x$ and $\mathcal{%
P}^{\ast }$ the vector space of all linear functionals on $\mathcal{P}$. If $%
L\in \mathcal{P}^{\ast }$, following Dirac we will denote the action of $L\,$%
on $p\left( x\right) \in \mathcal{P}\,\,$by $\left\langle L\right\vert
\left. p\left( x\right) \right\rangle $. A remarkable fact of the finite
operator calculus is that any element of $\mathcal{F\,}$can play a threefold
r\^{o}le: It can be regarded as a formal power series, as a linear
functional on $\mathcal{P}$ and also as a linear operator on $\mathcal{P}$.
To prove this, let us first notice that, given a formal power series
\begin{equation}
f\left( t\right) =\sum_{k=0}^{\infty }\frac{a_{k}}{k!}t^{k}\text{,}
\label{2.2}
\end{equation}%
we can associate it with a linear functional via the correspondence
\begin{equation}
\left\langle f\left( t\right) \right\vert \left. x^{n}\right\rangle =a_{n}%
\text{.}  \label{2.3}
\end{equation}%
Therefore
\begin{equation}
f\left( t\right) =\sum_{k=0}^{\infty }\frac{\left\langle f\left( t\right)
\right\vert \left. x^{k}\right\rangle }{k!}t^{k}\text{.}  \label{2.4}
\end{equation}%
In particular, since
\begin{equation}
\left\langle t^{k}\right\vert \left. x^{n}\right\rangle =n!\delta _{n,k}%
\text{,}  \label{2.5}
\end{equation}%
it follows that, for any polynomial $p\left( x\right) $,%
\begin{equation}
\left\langle t^{k}\right\vert \left. p\left( x\right) \right\rangle
=p^{\left( k\right) }\left( 0\right) \text{,}  \label{2.6}
\end{equation}%
where $p^{\left( k\right) }\left( 0\right) $ denotes the $k$--derivative of $%
p\left( x\right) $ evaluated at $x=0$. It is easily shown that any linear
functional $L\in \mathcal{P}^{\ast }\,$is of the form (\ref{2.4}).

If we interpret $t^{k}$ as the $k$--$th$ order derivative operator on $%
\mathcal{P}$, given a polynomial $p\left( x\right) $ we have that $%
t^{k}p\left( x\right) =p^{\left( k\right) }\left( x\right) $. Hence the
formal series
\begin{equation}
f\left( t\right) =\sum_{k=0}^{\infty }\frac{a_{k}t^{k}}{k!}  \label{3.6}
\end{equation}
is also regarded as a linear operator acting on $\mathcal{P}$. Depending on
the context, $t$ will play the r\^{o}le of a formal variable or that of a
derivative operator.

Now, some basic definitions and theorems of finite operator theory,
necessary in the subsequent considerations, are in order.

\begin{definition}
An operator $S\,$ commuting with the shift--operator $T$, i.e. $\left[ S,T%
\right] =0$, is said to be shift--invariant.
\end{definition}

Relevant examples of operators belonging to this class are provided by the
delta operators.

\begin{definition}
A delta operator $Q\,$is a shift--invariant operator such that $Qx=const\neq
0.$
\end{definition}

As has been proved in \cite{Rota}, there is an isomorphism between the ring
of formal power series in a variable $t$ and the ring of shift--invariant
operators, carrying
\begin{equation}
f\left( t\right) =\sum_{k=0}^{\infty }\frac{a_{k}t^{k}}{k!}\text{ into }%
\sum_{k=0}^{\infty }\frac{a_{k}Q^{k}}{k!}\text{.}  \label{3.7}
\end{equation}

In the following let us denote by $p_{n}\left( x\right) $, $n=0,1,2,\ldots,$ a
sequence of polynomials in $x$, where $p_{i}\left( x\right) $ is of order $i$
for all $i$.

\begin{definition}
A polynomial sequence $p_{n}\left( x\right) $ is called a sequence of basic
polynomials for a delta operator $Q$ if it satisfies the following
conditions:

%

1) $p_{0}\left( x\right) =1$;

2) $p_{n}\left( 0\right) =0$ \ for all $n>0$;

3) $Qp_{n}\left( x\right) =np_{n-1}\left( x\right) $.
\end{definition}

For each delta operator there exists a unique sequence of associated basic
polynomials.

\begin{definition}
A polynomial sequence $s_{n}\left( x\right) $ is called a set of Sheffer
polynomials for the delta operator $Q$ if

%

1) $s_{0}\left( x\right) =c\neq 0$;

2) $Qs_{n}\left( x\right) =ns_{n-1}\left( x\right) $.
\end{definition}

\begin{definition}
An Appell sequence of polynomials is a Sheffer sequence for the delta
operator $D$.
\end{definition}

Any shift invariant operator $S$ can be expanded into a formal power series
in terms of a delta operator $Q$:
\begin{equation*}
S=\sum_{k\geq 0}\frac{a_{k}}{k!}Q^{k}\text{,}
\end{equation*}
with $a_{k}=\left[ Sp_{k}\left( x\right) \right] \mid _{x=0}$, where $%
p_{k}\,\,$is the basic polynomial of order $k\,$ associated with $Q$. By
using the isomorphism (\ref{3.7}), a formal power series $s\left( t\right) $
is defined, which is called the\textit{\ indicator of }$S$.

\begin{remark}
A shift invariant operator is a delta operator if and only if it
corresponds, under the isomorphism (\ref{3.7}), to a formal power series $%
G\left( t\right) $ such that $G\left( 0\right) =0$ and $G^{\prime }\left(
0\right) \neq 0$. This series admits a unique compositional inverse.
\end{remark}

The umbral formalism also allows us to characterize the generating functions
of polynomial sequences of Sheffer and Appell type.

If $s_{n}\left( x\right) $ is a Sheffer sequence for the operator $Q$, then
there exists an invertible shift--invariant operator $S$ such that
\begin{equation*}
S\text{ }s_{n}\left( x\right) =q_{n}\left( x\right) \text{,}
\end{equation*}%
where $q_{n}\left( x\right) $ is a basic set for $Q$.

Let us denote by $s\left( t\right) $ and $q\left( t\right) $ the indicators
of the operators $S$ and $Q$. We will say that $s_{n}\left( x\right) $ is
the Sheffer sequence associated with $\left( s\left( t\right) ,q\left(
t\right) \right) $. The following result holds:
\begin{equation}
\frac{1}{s\left( q^{-1}\left( t\right) \right) }e^{xq^{-1}\left( t\right)
}=\sum_{n\geq 0}\frac{s_{n}\left( x\right) }{n!}t^{n}\text{,}  \label{2.19}
\end{equation}
where $q^{-1}\left( t\right) $ denotes the compositional inverse of $q\left(
t\right) $.

Another characterization of a Sheffer sequence is provided by the identity
\begin{equation*}
s_{n}\left( x+y\right) =\sum_{k\geq 0}\binom{n}{k} q_{k}\left( y\right)
s_{n-k}\left( x\right) \text{.}
\end{equation*}
In the specific case of the Appell sequences, which obey the equation
\begin{equation}
D\text{ }s_{n}\left( x\right) =ns_{n-1}\left( x\right) \text{,}  \label{2.21}
\end{equation}
there exists an invertible operator $g\left( t\right) $ such that
\begin{equation}
g\left( t\right) s_{n}\left( x\right) =x^{n}\text{.}  \label{2.22}
\end{equation}
Therefore, we will say that $s_{n}\left( x\right) $ is the Appell sequence
associated with $g\left( t\right) $. Its generating function is defined by

\begin{equation}
\sum_{k=0}^{\infty }\frac{s_{k}\left( x\right) }{k!}t^{k}=\frac{1}{g\left(
t\right) }e^{xt}\text{.}  \label{2.23}
\end{equation}
The Appell identity is
\begin{equation}
s_{n}\left( x+y\right) =\sum_{k=0}^{n}\binom{n}{k}s_{k}\left( y\right)
x^{n-k}\text{.}  \label{2.24}
\end{equation}

We also recall the polynomial expansion theorem. If $s_{n}\left( x\right) $
is a Sheffer set for the operators $\left( s\left( t\right) ,q\left(
t\right) \right) $, then for any polynomial $p\left( x\right) $ we have:
\begin{equation}
p\left( x\right) =\sum_{k\geq 0}\frac{\left\langle s\left( t\right) q\left(
t\right) ^{k}\mid p\left( x\right) \right\rangle }{k!}s_{k}\left( x\right)
\text{.}  \label{2.26}
\end{equation}%
In particular, if $q_{n}\left( x\right) $ is a basic sequence for $q(t)$, this formula
reduces to
\begin{equation}
p\left( x\right) =\sum_{k\geq 0}\frac{\left\langle q\left( t\right) ^{k}\mid
p\left( x\right) \right\rangle }{k!}q_{k}\left( x\right) \text{.}
\label{2.27}
\end{equation}

\section{Difference delta operators}

In this paper, we will consider a specific class of difference operators,
introduced in \cite{LTW1}, having the general form
\begin{equation}
Q\equiv \Delta _{p}=\frac{1}{\sigma }\sum_{k=l}^{m}a_{k}T^{k}\text{,}\qquad l%
\text{, }m\in \mathbb{Z}\text{,}\mathbb{\qquad }l<m\text{,\qquad }m-l=p\text{%
,}  \label{2.9}
\end{equation}%
where $a_{k}\,$and $\sigma $ are constants. In order to satisfy the
definition of delta operator, we must assume that
\begin{equation}
\sum_{k=l}^{m}a_{k}=0\text{,}\quad \sum_{k=l}^{m}ka_{k}=c\text{.}
\label{2.10}
\end{equation}%
We also require that in the continuous limit analogue, $\Delta _{p}\,$%
reproduces the standard derivative $D$; this implies $c=1$, i.e.%
\begin{equation}
\sum_{k=l}^{m}ka_{k}=1\text{.}  \label{2.11}
\end{equation}

\begin{definition}
A difference operator of the form (\ref{2.9}) which satisfies Eqs. (\ref%
{2.10}) and (\ref{2.11}) is called a delta operator of order $p=m-l$.
\end{definition}

We observe that Eq. (\ref{2.9}) involves $m-l+1\,$constants $a_{k}$, subject
to conditions (\ref{2.10}) and (\ref{2.11}). To fix all constants $%
a_{k}\,$we have to impose $m-l-1\,$further conditions. A possible choice is,
for instance,
\begin{equation}
\gamma _{p}:= \sum_{k=l}^{m}k^{p}a_{k}=0\text{,}\qquad p=2,3,...,m-l%
\text{.}  \label{2.12}
\end{equation}%
When conditions (\ref{2.10}), (\ref{2.11}) and (\ref{2.12}) are satisfied,
the operator (\ref{2.9}) provides an approximation of order $p$ of the
continuous derivative $D$, since
\begin{equation*}
{\Delta }_{p}{f}\mathop{\sim }_{\sigma \rightarrow 0}f^{\prime }\left(
x\right) +\frac{\sigma ^{m-l}}{\left( m-l+1\right) !}f^{\left( m-l-1\right)
}\left( x\right) \sum_{k=l}^{m}a_{k}k^{m-l-1}\text{.}
\end{equation*}%
From now on, we will put $\sigma =1$. The Pincherle derivative of a delta
operator is defined by the relation
\begin{equation*}
Q^{\prime }=\left[ Q,x\right] \text{.}
\end{equation*}%
Now, let us introduce a shift--invariant operator $\beta $ such that
\begin{equation}
\left[ \Delta _{n},x\beta \right] =1\text{.}  \label{2.13}
\end{equation}%
It follows that $\beta =\left( \Delta _{n}^{\prime }\right) ^{-1}$ (see \cite%
{LNO}--\cite{LTW2}). Such an operator is invertible \cite{Rota} and finite.
Indeed, using the identity $\beta \beta ^{-1}=1$, the action of $\beta $ on
a monomial of order $n$ is easily seen to be
\begin{equation}
\beta x^{n}=x^{n}-\sum_{j=0}^{n-1}\alpha _{j}^{n}x^{n-1-j}\text{,}
\label{2.14}
\end{equation}%
where $\alpha _{j}^{n}$ are defined via the recursion relation
\begin{equation*}
\alpha _{j}^{n}={\binom{n}{j+1}}\gamma _{j+2}-\sum_{l=0}^{j-1}{\binom{n}{l+1}%
}\gamma _{l+2}\alpha _{j-l-1}^{n-l-1}\text{,}
\end{equation*}%
with $\gamma _{j}=\sum_{k=l}^{m}k^{j}a_{k}$. We see that $\beta $ preserves
polynomial structures. When $\Delta =D,$ $\beta =1$. Other specific cases
are listed below: for
\begin{equation}
\Delta =\Delta ^{+}=T-1\text{,}\quad \beta =T^{-1}\text{;}
\end{equation}%
for
\begin{equation}
\Delta =\Delta ^{-}=1-T^{-1}\text{,}\quad \beta =T\text{.}
\end{equation}%
When
\begin{equation}
\Delta _{2}=\Delta ^{s}=\frac{T-T^{-1}}{2}\text{,}\quad \beta =\left( \frac{%
T+T^{-1}}{2}\right) ^{-1}\text{.}
\end{equation}%
Other nontrivial examples of higher--order operators are provided by
\begin{equation}
\Delta _{3}=-\left( T^{2}-2T+T^{-1}\right) ,  \label{2.30}
\end{equation}

\begin{equation}
\Delta _{4}=T^{2}-\frac{3}{2}T+\frac{3}{2}T^{-1}-T^{-2}\text{,}  \label{2.31}
\end{equation}

\begin{equation}
\Delta _{5}=T^{3}-2T^{2}+2T-2T^{-1}+T^{-2}\text{,}  \label{2.32}
\end{equation}

\begin{equation}
\Delta _{7}=T^{4}-T^{3}+T^{2}-2T+T^{-1}-T^{-2}+T^{-3}\text{.}  \label{2.34}
\end{equation}

Finally, we recall that $Q$ is a delta operator if and only if is of the
form $Q=DP$, where $P$ is an invertible, shift--invariant operator \cite{Rota}%
. It follows that the sequence of basic polynomials for $Q\,$ is expressible
in the form
\begin{equation}
p_{n}\left( x\right) =xP^{-n}x^{n-1}\text{.}  \label{2.35}
\end{equation}

\section{A generalization of the Stirling numbers}

According to Rota's approach, given any delta operator of the form (\ref{2.9}%
)--(\ref{2.11}) it is possible to associate  it with a sequence of basic polynomials.
In perfect analogy to the classical theory of Stirling numbers, this
allows us to study generalizations of Stirling numbers. In this section,
several concrete examples of basic polynomials as well as generalized
Stirling numbers are obtained. This construction represents a natural
realization of the very general scheme proposed by Ray in \cite{Ray}, who
first studied the basic sequences and defined the universal Stirling numbers
for any delta operator constructed via the complex homology theory. Here we
will follow a different but equivalent philosophy: the basic sequences for
the operators (\ref{2.9})--(\ref{2.11}) are simply derived from the
\textquotedblright discrete\textquotedblright\ representations of the
Heisenberg--Weyl algebra; the associated Stirling numbers furnish the
connection coefficients between the basic sequences and the standard power
sequence $\left\{ x^{n}\right\} _{n\in \mathbb{N}}$.

We recall that the classical Stirling numbers of the first and second kind $%
s\left( n,k\right) \,$ and $S\left( n,k\right) $ are defined by the
relations (see, e.g., \cite{Riordan1})
\begin{equation}
\left( x\right) _{n}=\sum_{k=0}^{n}s\left( n,k\right) x^{k}  \label{3.1}
\end{equation}%
and
\begin{equation}
x^{n}=\sum_{k=0}^{n}S\left( n,k\right) \left( x\right) _{k}\text{,}
\label{3.2}
\end{equation}%
respectively. Here $\left( x\right) _{n}=x\left( x-1\right) \ldots\left(
x-n+1\right) $ denotes the lower factorial polynomial of order $n$. The
numbers $S\left( n,k\right) \,$also admit a representation in terms of a
generating function:
\begin{equation*}
\sum_{n=k}^{\infty }S\left( n,k\right) \frac{t^{n}}{n!}=\frac{\left(
e^{t}-1\right) ^{k}}{k!}\text{.}
\end{equation*}%
As is well known, the lower factorial polynomials are the basic sequence
associated with the forward derivative: $\Delta \left( x\right) _{n}=n\left(
x\right) _{n-1}$. More generally, the basic sequence associated with the
operator $\Delta _{p}\,$, denoted by $(x)_{n}^{p}$, is \cite{LTW1}
\begin{equation}
(x)_{n}^{p}=\left( x\beta \right) ^{n}1\text{.}  \label{3.3}
\end{equation}%
Therefore,
\begin{equation}
\Delta _{p}(x)_{n}^{p}=n(x)_{n-1}^{p}\text{.}  \label{3.3a}
\end{equation}%
Explicitly, the first basic polynomials are expressed by
\begin{equation*}
(x)_{0}^{p}=\left( x\beta \right) ^{0}1=1\text{,}
\end{equation*}%
\begin{equation*}
(x)_{1}^{p}=\left( x\beta \right) ^{1}1=x\text{,}
\end{equation*}%
\begin{equation*}
(x)_{2}^{p}=\left( x\beta \right) ^{2}1=x^{2}-\gamma _{2}x\text{,}
\end{equation*}%
\begin{equation*}
(x)_{3}^{p}=\left( x\beta \right) ^{3}1=x^{3}-3\gamma _{2}x^{2}-\left(
\gamma _{3}-3\gamma _{2}^{2}\right) x\text{,}
\end{equation*}%
\begin{equation}
(x)_{4}^{p}=\left( x\beta \right) ^{4}1=x^{4}-6\gamma _{2}x^{3}+\left(
-4\gamma _{3}+15\gamma _{2}^{2}\right) x^{2}+\left( -\gamma _{4}+10\gamma
_{2}\gamma _{3}-15\gamma _{2}^{3}\right) x\text{,}
\end{equation}%
and so on, with $\gamma _{j}$ defined by
\begin{equation}
\gamma _{j}=\sum_{k=l}^{m}a_{k}k^{j},\quad \gamma _{0}=0,\quad \gamma
_{1}=1,\quad j=0,1,2,\ldots\text{.}
\end{equation}

\begin{definition}
The \textit{generalized Stirling numbers of the first kind} and order $p$
associated with the operators (\ref{2.9})--(\ref{2.11}), denoted by $%
s^{p}\left( n,k\right) $, are defined by
\begin{equation}
(x)_{n}^{p}=\sum_{k=0}^{n}s^{p}\left( n,k\right) x^{k}\text{.}  \label{GS1}
\end{equation}
\end{definition}

\begin{remark}
In the subsequent considerations, with an abuse of notation, we will use
Roman numerals for denoting the values of $p$, when $p$ appears as an index
in sequences of polynomials or numbers. 
\end{remark}

Observe that $(x)_{n}^{I}=\left( x\right) _{n}$ and $s^{I}\left( n,k\right)
=s\left( n,k\right) $. Since for any polynomial $p\left( x\right) $ the
expansion
\begin{equation}
p\left( x\right) =\sum_{k\geqslant 0}\frac{\left\langle t^{k}\mid p\left(
x\right) \right\rangle }{k!}x^{k}\text{,}
\end{equation}%
holds, it follows that
\begin{equation}
s^{p}\left( n,k\right) =\frac{1}{k!}\left\langle t^{k}\mid
(x)_{n}^{p}\right\rangle \text{.}
\end{equation}

\begin{definition}
The \textit{generalized Stirling numbers of the second kind and order p},
denoted by $S^{p}\left( n,k\right) $, are defined by
\begin{equation}
x^{n}=\sum_{k=0}^{n}S^{p}\left( n,k\right) (x)_{k}^{p}\text{.}  \label{GS2}
\end{equation}
\end{definition}

They admit the generating function
\begin{equation}
\sum_{n=k}^{\infty }S^{p}\left( n,k\right) \frac{t^{n}}{n!}=\frac{\left[
\Delta _{p}\left( t\right) \right] ^{k}}{k!}\text{,}
\end{equation}%
where $\Delta _{p}\left( t\right) =\sum_{j}a_{j}e^{jt}\,$is the indicator of
the operator $\Delta _{p}$. As a consequence of formula (\ref{2.27}), one
immediately obtains that
\begin{equation}
S^{p}\left( n,k\right) =\frac{1}{k!}\left\langle \left(
\sum_{j=l}^{m}a_{j}e^{jt}\right) ^{k}\mid x^{n}\right\rangle \text{,}\quad
\text{ }m-l=p\text{.}
\end{equation}%
The polynomials $(x)_{n}^{p}\,$satisfy the binomial identity
\begin{equation}
\left( x+y\right) _{n}^{p}=\sum_{k=0}^{n}\binom{n}{k}(x)_{k}^{p}%
\,(y)_{n-k}^{p}\text{.}
\end{equation}%
A relevant feature of relations (\ref{GS1}) and (\ref{GS2}) is that they are
inverse to each other. This immediately implies that
\begin{equation}
\delta _{mn}=\sum_{k}s^{p}\left( n,k\right) S^{p}\left( k,m\right)
=\sum_{k}S^{p}\left( n,k\right) s^{p}\left( k,m\right) \text{.}  \label{3.10}
\end{equation}

\begin{remark}
Due to formulae (\ref{2.13}) and (\ref{3.3a}), it is natural to interpret
the operators $\Delta $ and $x\beta $ as quantum mechanical annihilation and
creation operators $\ a$ and $\ a^{\dagger }$ acting on a finitely generated
space, as noticed in \cite{T1}, \cite{T2}. For instance,
\begin{equation*}
\left( x\beta \right) ^{n}\cdot 1=\left( a^{\dagger }\right) ^{n}\left\vert
0\right\rangle =\left\vert n\right\rangle \text{.}
\end{equation*}%
\begin{equation*}
\Delta _{p}(x\beta )^{n}\cdot 1=n\left\vert n-1\right\rangle \text{.}
\end{equation*}
\end{remark}

\section{Bernoulli--type polynomials and numbers of the first kind}

In this section, some realizations of the polynomial class (\ref{I.7}) are
presented. We shall restrict to sequences of polynomials obtainable in the
umbral context described above. We propose the following definition.

\begin{definition}
The higher-order \textit{Bernoulli--type polynomials of the first kind} are
the polynomials generated by means of the relation
\begin{equation}
\sum_{k=0}^{\infty }\frac{\,B_{k,a}^{p}\left( x\right) }{k!}%
t^{k}=(J_{p}\left( t\right))^{a}e^{xt}\text{,} 
\label{4.1}
\end{equation}%
with $a\neq0$ and
\begin{equation}
J_{p}\left( t\right) =\frac{t}{\sum_{j=l}^{m}a_{j}e^{tj}}\text{,}
\label{4.2}
\end{equation}%
where $\sum_{j=l}^{m}a_{j}e^{tj}\,$is the indicator of a delta operator of
order $p$.
\end{definition}

\begin{remark}
Usually, we will assume that $a_{j}\in \mathbb{Q}$, as turns out to be the case in many applications. However, more
general operator structures with real coefficients $a_{j}$ could be
considered within the same framework.

In the following, when $a=1$, the index $a$ will be omitted for simplicity.
\end{remark}
\begin{remark}
When $p=a=1$, $J_1(t)=\frac{t}{e^{t}-1}$, and in this case we obtain
the classical Bernoulli polynomials. In the case $a\neq 1$, $p=1$, the
polynomials (\ref{4.2}) reduce to the standard higher--order Bernoulli
polynomials. 
\end{remark}
\begin{definition}
\smallskip The Bernoulli--type numbers $B_{k,a}^{p}$ are defined by the relation
\begin{equation}
\sum_{k=0}^{\infty }\frac{\,B_{k,a}^{p}}{k!}t^{k}=J_{p}\left( t\right) ^{a}, \qquad a\neq 0%
\text{.}  \label{4.3}
\end{equation}
\end{definition}

\begin{remark}
The operators $J_{p}\left( t\right) \,$ and the corresponding numbers for $%
a=1$ have been studied in \cite{Ray} in the context of complex oriented homology theory.
Here again the construction of these operators is realized using a
difference operator approach.
\end{remark}

\begin{remark}
Let us introduce the operator $\Delta^{a}_p$, which is defined to be the operator whose indicator is $\Delta^{a}_{p}(t):=t (J_p(t))^{-a}$. In the special case $a=1$,  we simply get the difference operator $\Delta^{1}_p:=\Delta_p$. 
\end{remark}
\smallskip We can now prove a result stated in the Introduction 
\footnote{This result has been extended now to the general case of the higher-order polynomials $B^p_{m,a}(x)$. 
In the published version,  the symbol $B^p_{m,a}(x)$ was used (due to an unfortunate notational misprint) instead of $B^{p}_{m}(x)$.  The present notation is the correct one.}

.

\begin{lemma}
\smallskip The polynomials $B_{m,a}^{p}\left( x\right) $ are 
determined by the two properties
\begin{equation}
DB_{m,a}^{p}\left( x\right) =mB_{m-1,a}^{p}\left( x\right) \text{,}
\label{4.12a}
\end{equation}%
\begin{equation}
\Delta^{a}_{p}B_{m,a}^{p}\left( x\right) =mx^{m-1}\text{,}  \label{4.12b}
\end{equation}
with the condition $B_{0,a}^{p}\left( x\right) =1$. 
\end{lemma}

\begin{proof} The argument we propose is very simple and closer to the spirit of the finite operator theory \cite{Roman}. Identity (\ref{4.12a}) follows immediately from (\ref{4.1}).
 Observe that, being an Appell sequence, the higher-order Bernoulli-type polynomials satisfy the relation \beq
 B^{p}_{m,a}(x)= (J_p(t))^{a} x^m.
\eeq  The result follows by applying to both sides the operator $\Delta^{a}_{p}(t)$, taking into account the initial condition. 
\end{proof}

\begin{remark}
For $a=1$, we obtain a new family of linear difference equations involving the Bernoulli-type polynomials. Also, for $a=1$, $p=1$, $\Delta^{1}_1\equiv \Delta$ so that we obtain again the conditions characterizing the classical Bernoulli polynomials. 
\end{remark}
\textbf{Examples}. The case $p=2$. We get a sequence of polynomials which we
will call the \textit{central Bernoulli polynomials} (of the first kind).
The generating function reads
\begin{equation}
\sum_{k=0}^{\infty }\frac{B_{k,a}^{II}\left( x\right) }{k!}t^{k}=\left(
\frac{t}{\sinh t}\right) ^{a}e^{xt}\text{.}  \label{4.11}
\end{equation}

For $x=0$ and $a=1$, one obtains the numbers $B_{n}^{II}$. From the generating function
it emerges that $B_{k}^{II}=0$ for $k$ odd. The first central Bernoulli
numbers $B_{2k}^{II}$ are: $1,-\frac{1}{3},\frac{7}{15},-\frac{31}{21},\frac{%
127}{15},-\frac{2555}{33},$ etc.

The polynomials (\ref{4.11}) can be directly expressed in terms of the
classical Bernoulli polynomials:
\begin{equation*}
B_{n}^{II}\left( x\right) =2^{n}B_{n}^{I}\left( \frac{x+1}{2}\right) \text{,}%
\quad B_{n}^{I}\left( x\right) =\frac{B_{n}^{II}\left( 2x-1\right) }{2^{n}}%
\text{.}
\end{equation*}

Other properties are now briefly discussed. From the polynomial expansion
theorem we get
\begin{equation}
B_{n}^{II}\left( x\right) =\sum_{k=0}^{n}\frac{\left\langle \left( \frac{%
e^{t}-e^{-t}}{2}\right) ^{k}\mid B_{n,a}^{II}\left( x\right) \right\rangle }{%
k!}\left( x\right) _{k}^{II}\text{.}
\end{equation}%
Since
\begin{equation*}
\left\langle \left( \frac{e^{t}-e^{-t}}{2}\right)^{k}\mid
B_{n,1}^{II}\left( x\right) \right\rangle =\left\langle \left( \frac{%
e^{t}-e^{-t}}{2}\right) ^{k}\mid \frac{2t}{e^{t}-e^{-t}}x^{n}\right\rangle =
\end{equation*}%
\begin{equation}
\left\langle \left( \frac{e^{t}-e^{-t}}{2}\right) ^{k-1}\mid
nx^{n-1}\right\rangle =n\left( k-1\right) !S^{II}\left( n-1,k-1\right) \text{%
,}
\end{equation}%
we have:
\begin{equation}
B_{n}^{II}\left( x\right) =B_{n}^{II}\left( 0\right) +\sum_{k=1}^{n}%
\frac{n}{k}S^{II}\left( n-1,k-1\right)(x)_{k}^{II} \text{.}  \label{6.15}
\end{equation}%
For $p=3$, the corresponding delta operator is given by formula (\ref{2.30}%
), with the generating function
\begin{equation}
\sum_{k=0}^{\infty }\frac{B_{k}^{III}\left( x\right) }{k!}t^{k}=-\frac{%
te^{xt}}{e^{2t}-2e^{t}+e^{-t}}\text{.}  \label{BIII}
\end{equation}%
For the other cases, we get similar results. For instance,
\begin{equation}
\sum_{k=0}^{\infty }\frac{B_{k}^{V}\left( x\right) }{k!}t^{k}=\frac{te^{xt}}{%
e^{3t}-2e^{2t}+2e^{t}-2e^{-t}+e^{-2t}}\text{,}  \label{BV}
\end{equation}%
\begin{equation}
\sum_{k=0}^{\infty }\frac{B_{k}^{VII}\left( x\right) }{k!}t^{k}=\frac{%
-te^{xt}}{e^{4t}-e^{3t}+e^{2t}-2e^{t}+e^{-t}-e^{-2t}+e^{-3t}}\text{,}
\label{BVII}
\end{equation}%
and so on. As an example, the first polynomials $B_{n}^{III}\left( x\right)
\,$ are given by
\begin{equation*}
B_{0}^{III}\left( x\right) =1\text{, }B_{1}^{III}\left( x\right) =x+\frac{3}{%
2}\text{, }B_{2}^{III}\left( x\right) =x^{2}+3x+\frac{37}{6}\text{,}
\end{equation*}%
\begin{equation*}
B_{3}^{III}\left( x\right) =x^{3}+\frac{9}{2}x^{2}+\frac{37}{2}x+39\text{,}
\end{equation*}%
\begin{equation*}
B_{4}^{III}\left( x\right) =x^{4}+6x^{3}+37x^{2}+156x+\frac{9719}{30}\text{,}
\end{equation*}%
\begin{equation*}
B_{5}^{III}\left( x\right) =x^{5}+\frac{15}{2}x^{4}+\frac{185}{3}%
x^{3}+390x^{2}+\frac{9719}{6}x+3365\text{,}
\end{equation*}%
\begin{equation*}
B_{6}^{III}\left( x\right) =x^{6}+9x^{5}+\frac{185}{2}x^{4}+780x^{3}+\frac{%
9719}{2}x^{2}+20190x+\frac{1762237}{42}\text{,}\ldots
\end{equation*}


The relation (\ref{6.15}) among higher order Stirling--type numbers,
Bernoulli--type numbers and Bernoulli--type polynomials can be generalized
as follows:
\begin{equation}
B_{n}^{p}\left( x\right) =B_{n}^{p}\left( 0\right) +\sum_{k=1}^{n}\frac{n%
}{k}S^{p}\left( n-1,k-1\right)(x_k)^{p} \text{.}
\end{equation}%
The recurrence relation for the polynomials $B_{n,a}^{p}\left( x\right) $ is
formally derived from the defining relation (\ref{4.1}):
\begin{equation}
\sum_{j=l}^{m}a_{j}e^{jt}B_{n,a}^{p}\left( x\right) =\left( \frac{t}{%
\sum_{j=l}^{m}a_{j}e^{tj}}\right) ^{a-1}tx^{n}=nB_{n-1,a-1}^{p}\left(
x\right) \text{.}
\end{equation}%
Therefore, the difference equation solved by the Bernoulli--type polynomials
is, for every $p$,%
\begin{equation}
\sum_{j=l}^{m}a_{j}B_{n,a}^{p}\left( x+j\right) =nB_{n-1,a-1}^{p}\left(
x\right) \text{,}
\end{equation}%
where the $a_{k}$ satisfy the constraints (\ref{2.10}), (\ref{2.11})
. When $p=2$, we obtain
\begin{equation} \label{eq.:5.16}
B_{n,a}^{II}\left( x+1\right) -B_{n,a}^{II}\left( x-1\right)
-2 n B_{n-1,a-1}^{II}\left( x\right) =0\text{.}
\end{equation}%
Also we have the recurrence
\begin{equation}
B_{n+1,a}^{p}\left( x\right) =\left( x-\frac{g^{\prime }\left( t\right) }{%
g\left( t\right) }\right) B_{n,a}^{p}\left( x\right) \text{,}
\end{equation}%
where $g(t)$ is the operator
\begin{equation*}
g\left( t\right) =\left( \frac{\sum_{k=l}^{m}a_{k}e^{kt}}{t}\right) ^{a}%
\text{.}
\end{equation*}%
Using the relations $\left[ t,x\right] =1$ and $B_{n,a}^{p}(x)=g\left( t\right)
^{-1}x^{n}$, which is a consequence of Eq. (\ref{2.22}), we get
\begin{equation}
\left( n+1\right) B_{n,a}^{p}\left( x\right) =\left( xt+1-a\frac{%
\sum_{k=l}^{m}a_{k}e^{kt}\left( tk-1\right) }{\sum_{k=l}^{m}a_{k}e^{kt}}%
\right) B_{n,a}^{p}\left( x\right) \text{.}
\end{equation}%
From the previous equation we derive the formula expressing the
Bernoulli--type polynomials of order $a+1\,$in terms of those of order $a$:
\begin{equation}
\sum_{k}ka_{k}B_{n,a+1}^{p}\left( x+k\right) =\left( 1-\frac{n}{a}\right)
B_{n,a}^{p}\left( x\right) +\frac{nx}{a}B_{n-1,a}^{p}\left( x\right) \text{.}
\end{equation}%
Observe that, in accordance with Eq. (\ref{2.35}), the polynomials
\begin{equation}
p_{n}\left( x\right) =x\left( \frac{t}{\sum_{k}a_{k}e^{kt}}\right)
^{na}x^{n-1}=xB_{n-1,\text{ }na}^{p}\left( x\right)
\end{equation}%
represent the basic sequence associated with the operator $Q=DP=D \left( \frac{%
\Delta _{p}}{D}\right) ^{a}\equiv\Delta_p^{a}$. Consequently, for $a=1$, we deduce the
relation expressing the basic sequence $(x)_{n}^{p}$ in terms of
Bernoulli--type polynomials:
\begin{equation} \label{eq:5.22}
(x)_{n}^{p}=xB_{n-1,n}^{p}\left( x\right) \text{.}
\end{equation}

An interesting connection between the Bernoulli numbers and the Stirling
numbers of the same order holds. Indeed, using Eq. \eqref{eq:5.22} and its differential consequences, after some algebraic manipulations we get
\begin{equation*}
s^{p}\left( n,r\right) =\frac{1}{r!}\left\langle t^{r}\mid
(x)_{n}^{p}\right\rangle =\frac{1}{r!}\left\langle t^{0}\mid
t^{r}(x)_{n}^{p}\right\rangle
\end{equation*}%
\begin{equation}
=\binom{n-1}{r-1}B_{n-r,n}^{p}\left( 0\right) \text{.}
\end{equation}%
Therefore,
\beq \label{eq:5.24}
(x)^{p}_{n}= \sum_{r=0}^n \binom{n-1}{r-1}B^{p}_{n-r,n}(0)x^r \ .
\eeq
Analogously,
\begin{equation}
S^{p}\left( n,r\right) =\frac{1}{r!}\left\langle \left( \frac{%
\sum_{k=l}^{m}a_{k}e^{kt}}{t}\right) ^{r}\mid t^{r}x^{n}\right\rangle =%
\binom{n}{r}B_{n-r,-r}^{p}\left( 0\right) \text{,}
\end{equation}%
and
\begin{equation}
x^{n}=\sum_{r=0}^{n}\binom{n}{r}B_{n-r,-r}^{p}\left( 0\right) (x)_{r}^{p}%
\text{.}
\end{equation}

\section{New Bernoulli--type polynomials of the second kind}

By analogy with the standard Bernoulli polynomials of the second kind (see, e.g., \cite{Roman}), we introduce the following sequences of polynomials.

\begin{definition}
The higher--order Bernoulli--type polynomials of the second kind are the polynomials defined by
\begin{equation}
b_{n}^{p}\left( x\right) =J_{p}\left( t\right)^{-1} (x)_{n}^{p},\quad n\in
\mathbb{N} \text{,}  \label{5.1}
\end{equation}%
where $J_{p}\left( t\right) $ is the operator  (\ref{4.2}).
\end{definition}
These polynomials represent a Sheffer sequence for the operator $\Delta _{p}$ of Eqs. (\ref{2.9}), (\ref{2.10}):
\begin{equation}
\Delta _{p}\,b_{n}^{p}\left( x\right) =nb_{n-1}^{p}\left( x\right) \text{,}
\label{5.2}
\end{equation}%
and therefore they satisfy the identity
\begin{equation}
b_{n}^{p}\left( x+y\right) =\sum_{k=0}^{n}\binom{n}{k}b_{k}^{p}\left(
y\right) (x)_{n-k}^{p}
\end{equation}%
which relates them with the higher factorial polynomials. In particular, for
$y=0$, we get
\begin{equation}
b_{n}^{p}\left( x\right) =\sum_{k=0}^{n}\binom{n}{k}b_{k}^{p}\left( 0\right)
(x)_{n-k}^{p}\text{.}
\end{equation}%
From Eq. (\ref{5.2}) we deduce the difference equation satisfied by the
polynomials (\ref{5.1}):
\begin{equation}
\sum_{k}a_{k}b_{n}^{p}\left( x+k\right) =nb_{n-1}^{p}\left( x\right) \text{.}
\label{5.5}
\end{equation}%
The case $J_1(t)^{-1}=\frac{e^{t}-1}{t}$ just reproduces the standard Bernoulli polynomials of the
second kind. The generating function of the sequences \eqref{5.1} can be obtained in some specific cases.
Let us briefly discuss the second-order case. Since the operator $\mathfrak{J}^{II}:=%
\frac{e^{t}-e^{-t}}{2t}$ acts as follows:
\begin{equation*}
\mathfrak{J}^{II}p\left( x\right) =\frac{1}{2}\int_{x-1}^{x+1}p\left( u\right) du\text{,%
}
\end{equation*}%
we deduce an explicit expression for the central Bernoulli polynomials of
the second kind:
\begin{equation}
b_{n}^{II}\left( x\right) =\mathfrak{J}^{II}\left( x\right) _{n}^{II}=\frac{1}{2}%
\int_{x-1}^{x+1}(u)_{n}^{II}du\text{.}  \label{bII}
\end{equation}%
The generating function is
\begin{equation}
\sum_{k=0}^{\infty }\frac{b_{k}^{II}\left( x\right) t^{k}}{k!}=\frac{t}{\log
\left( t+\sqrt{1+t^{2}}\right) }\left( t+\sqrt{1+t^{2}}\right) ^{x}\text{.}
\end{equation}%
For $x=0$ we get%
\begin{equation}
\sum_{k=0}^{\infty }\frac{b_{k}^{II}\left( 0\right) t^{k}}{k!}=\frac{t}{\log
\left( t+\sqrt{1+t^{2}}\right) }\text{.}
\end{equation}%
Let us derive the recurrence relation satisfied by the polynomials $%
b_{n}^{II}\left( x\right) $. From Eq. (\ref{bII}) we get

\begin{equation}
tb_{n}^{II}\left( x\right) =\left( \frac{e^{t}-e^{-t}}{2}\right)
(x)_{n}^{II}=n(x)_{n-1}^{II}\text{.}
\end{equation}
Integrating we obtain
\begin{equation}
b_{n}^{II}\left( x\right) -b_{n}^{II}\left( 0\right)
=n\int_{0}^{x}(u)_{n-1}^{II}du\text{.}
\end{equation}
\smallskip The difference equation (\ref{5.5}) reduces to

\begin{equation}
b_{n}^{II}\left( x+1\right) =b_{n}^{II}\left( x-1\right)
+2nb_{n-1}^{II}\left( x\right) \text{.}
\end{equation}%
Since 
\begin{equation}
b_{n}^{II}\left( 1\right) =\frac{1}{2}\int_{0}^{2}(u)_{n-1}^{II}du\text{,}
\end{equation}%
another representation of $b^{II}_{n}(0)$ is
\begin{equation}
b_{n}^{II}\left( 0\right) =\frac{1}{2}\int_{0}^{2}(u)_{n-1}^{II}du-n\int_{0}^{1}\left(
u\right) _{n-1}^{II}du\text{.}
\end{equation}%
Notice the analogy with the classical Bernoulli numbers of the second kind,
defined by (see \cite{Roman})
\begin{equation}
b_{n}\left( 0\right) =\left\langle \frac{e^{t}-1}{t}\mid \left( x\right)
_{n}\right\rangle =\int_{0}^{1}\left( u\right) _{n}du\text{.}
\end{equation}%
There is a connection between Bernoulli--type polynomials of the second kind
and generalized Stirling numbers of the first kind. Since
\begin{equation*}
b_{n}^{p}\left( x\right) =\sum_{k=0}^{n}\frac{1}{k!}\left\langle t^{k}\mid
b_{n}^{p}\left( x\right) \right\rangle x^{k}\text{,}
\end{equation*}%
we get
\begin{equation}
b_{n}^{p}\left( x\right) =b_{n}^{p}\left( 0\right) +\sum_{k=1}^{n}\frac{n}{k}%
s^{p}\left( n-1,k-1\right) x^{k}\text{.}
\end{equation}

Many other properties and identities, which we will not discuss here for the sake of brevity, can be
derived using operator techniques.

\section{New Euler--type polynomials}

In this Section, a new class of polynomial sequence of Appell type is
discussed. For many aspects, it can be considered to be a natural
generalization of the Euler polynomials. The definition proposed in this
paper is different from that in \cite{You}.

Euler polynomials and numbers (introduced by Euler in 1740) also possess an
extensive literature and several interesting applications in Number Theory
(see, for instance, \cite{DSS}, \cite{Roman}, \cite{You}). In many respects,
they are closely related to the theory of Bernoulli polynomials and numbers.

\begin{definition}
The Euler--type polynomials are the Appell sequence generated by
\begin{equation}
\sum_{k=0}^{\infty }\frac{E_{k}^{p,a,\omega }\left( x\right) }{k!}%
t^{k}=\left( 1+\frac{\Delta_{p}(t)}{\omega }\right) ^{-a}e^{xt}\text{,}
\label{9.1}
\end{equation}%
where  $a,\omega\neq 0$ and $\Delta_p(t)$ is the indicator of the operator $\Delta_p$.
\end{definition}

For $a=1$, $\Delta_p(t)=e^{t}-1$, $\omega =2$, we obtain the classical Euler polynomials. It will
be shown that these new polynomials possess many of the properties of their
classical analogues.

We recall that the Euler polynomials of order $a$, which will be denoted by
the symbol $E_{n}^{a}\left( x\right) $, are the Appell sequence associated
with the operator
\begin{equation*}
g\left( t\right) =\left( 1+\frac{e^{t}-1}{2}\right) ^{a},\quad a\neq 0\text{.%
}
\end{equation*}%
Consequently
\begin{equation}
E_{n}^{a}\left( x\right) =\left( 1+\frac{e^{t}-1}{2}\right) ^{-a}x^{n}\text{.%
}
\end{equation}%
By analogy, we have the following

\begin{lemma}
The Euler--type polynomial sequences are given by the relation
\begin{equation}
E_{n}^{p,a,\omega }\left( x\right) =\left( 1+\frac{\Delta _{p}(t)}{\omega }%
\right) ^{-a}x^{n}\text{,}  \label{6.3}
\end{equation}%
where $a,\omega $ $\in \mathbb{R}$, $a,\omega\neq 0$.
\end{lemma}

For the sake of clarity, we will omit the superscript $\omega $. The
Euler--type polynomials satisfy the Appell property
\begin{equation*}
DE_{n}^{p,a}\left( x\right) =nE_{n-1}^{p,a}\left( x\right)
\end{equation*}%
and the Appell binomial identity

\begin{equation}
E_{n}^{p,a}\left( x+y\right) =\sum_{k=0}^{n}\binom{n}{k}E_{k}^{p,a}\left(
y\right) x^{n-k}\text{.}
\end{equation}%
For every $p$, the polynomials (\ref{9.1}) are an Appell
cross--sequence:
\begin{equation}
E_{n}^{p,\text{ }a+b}\left( x+y\right) =\sum_{k=0}^{n}\binom{n}{k}E_{k}^{p,%
\text{ }a}\left( y\right) E_{n-k}^{p,\text{ }b}\left( x\right) \text{.}
\end{equation}%
We also have a recurrence formula :
\begin{equation}
E_{n+1}^{p,a}\left( x\right) =\left( x-\frac{g_{p}^{\prime }\left( t\right)
}{g_{p}\left( t\right) }\right) E_{n}^{p,a}\left( x\right) \text{,}
\label{9.8}
\end{equation}%
where $g_{p}\left( t\right) =\left( 1+\frac{\Delta _{p}\left( t\right) }{%
\omega }\right)^{a} $.

There exists an interesting analogue of the Boole summation formula. Since
for any function $h\left( t\right) $,%
\begin{equation*}
h\left( t\right) =\sum_{k=0}^{\infty }\frac{\left\langle h\left( t\right)
\right\vert \left. E_{k}^{p,a}\left( x\right) \right\rangle }{k!}\left( 1+%
\frac{\Delta _{p}(t)}{\omega }\right) ^{a}t^{k}\text{,}
\end{equation*}%
then, choosing $h\left( t\right) =\exp \left( yt\right) $, we get
\begin{equation}
\exp \left( yt\right) =\sum_{k=0}^{\infty }\frac{E_{k}^{p,a}\left( y\right)
}{k!}\left( 1+\frac{\Delta _{p}(t)}{\omega }\right) ^{a}t^{k}\text{.}
\label{6.9}
\end{equation}%
Once applied to a polynomial $p\left( x\right) $, we obtain for any $p\,$
the expansion
\begin{equation*}
p\left( x+y\right) =\sum_{k=0}^{\infty }\frac{E_{k}^{p,a}\left( y\right) }{k!%
}\left( 1+\frac{\Delta _{p}(t)}{\omega }\right) ^{a}p^{\left( k\right) }\left(
0\right) \text{.}
\end{equation*}%
Let us consider in more detail the operator (for $p=2$)
\begin{equation}
g^{II}\left( t\right) =\left( 1+\frac{e^{t}-e^{-t}}{\omega }\right) ^{a}%
\text{.}
\end{equation}%
Correspondingly,
\begin{equation}
E_{n}^{II,a}\left( x\right) =\left( 1+\frac{e^{t}-e^{-t}}{\omega }\right)
^{-a}x^{n}\text{.}  \label{EII}
\end{equation}%
Since
\begin{equation}
g^{II}\left( t\right) ^{-1}=\left( \frac{1}{1+\frac{e^{t}-e^{-t}}{\omega }}%
\right) ^{a}=\sum_{j=0}^{\infty }{\binom{-a}{j}}\left( \frac{e^{t}-e^{-t}}{%
\omega }\right) ^{j}\text{,}
\end{equation}%
the polynomials (\ref{EII}) are of the form
\begin{equation}
E_{n}^{II,a}\left( x\right) =\sum_{j=0}^{\infty }{\binom{-a}{j}}\left( \frac{%
e^{t}-e^{-t}}{\omega }\right) ^{j}x^{n}=\sum_{k=0}^{n}\sum_{j=0}^{n}{\binom{%
-a}{j}}\frac{j!}{k!}\frac{2^{j}}{\omega ^{j}}S^{II}\left( k,j\right) x^{n-k}%
\text{.}  \label{9.12}
\end{equation}%
From this equation and from the definition (\ref{EII}) we derive the
difference equation
\begin{equation}
E_{n}^{II,a}\left( x+1\right) -E_{n}^{II,a}\left( x-1\right) =\omega
(E_{n}^{II,a-1}\left( x\right)-E_{n}^{II,a}\left(x\right)) \text{.}
\end{equation}%
The analogue of the Newton expansion is given by the formula
\begin{equation}
E_{n}^{II,a}\left( x\right) =\sum_{k=0}^{n}\sum_{j=k}^{n}{\binom{-a}{k}}%
\frac{2^{k}}{\omega ^{k}}\left( j\right) _{k}S^{II}\left( n,j\right) \left(
x\right) _{j-k}^{II}\text{,}  \label{9.14}
\end{equation}%
which is easily derived from eq. (\ref{9.12}) and from the formula
\begin{equation}
\left( \frac{e^{t}-e^{-t}}{2}\right) ^{k}x^{n}=\sum_{j=k}^{n}S^{II}\left(
n,k\right) \left( j\right) _{k}\left( x\right) _{j-k}^{II}\text{.}
\end{equation}%
Finally, from the recurrence formula (\ref{9.8}) (for $p=2$) we have:
\begin{equation*}
E_{n+1}^{II,a}\left( x\right) =\left( x-a\frac{e^{t}+e^{-t}}{\omega
+e^{t}-e^{-t}}\right) E_{n}^{II,a}\left( x\right) \text{,}
\end{equation*}%
providing another difference equation satisfied by $E_{n}^{II,a}\left(
x\right) $:
\begin{equation} \label{eq:7.15}
E_{n+1}^{II,a}\left( x\right) =xE_{n}^{II,a}\left( x\right) -\frac{a}{\omega
}\left( E_{n}^{II,a+1}\left( x+1\right) +E_{n}^{II,a+1}\left( x-1\right)
\right) \text{.}
\end{equation}
Other cases can be studied similarly.

\vspace{2mm}

\textbf{Acknowledgments}

\vspace{2mm}

I am grateful to A. Granville, who pointed out to me References \cite{Clarke}
and \cite{Ray} (after an early version of the paper was completed) and to B.
Dubrovin, R. A. Leo, D. Levi and P. Winternitz for many useful discussions.

\vspace{2mm}

\textbf{Errata corrige}

\vspace{2mm}

The following main changes have been done and typos corrected w.r.t. the published version:

\vspace{2mm}

(1) A new, more general statement and proof of Lemma 1 has been provided. In particular, in Eq. \eqref{4.12b} of Lemma 1, $\Delta_p \longrightarrow \Delta^{a}_p$.   

(2) The Remarks $5$, $7$ and $8$ have been added.

(3) The factor $2$ in the r.h.s. of Eq. \eqref{eq.:5.16} has been restored. 
The factor $a$ missing in the r.h.s. of Eq. \eqref{eq:7.15} has been restored. The r.h.s. of Eq. \eqref{eq:5.24} has been corrected.


\begin{thebibliography}{99}
\bibitem{Adel} A. Adelberg, Universal higher order Bernoulli numbers and
Kummer and related congruences, J. Number Theory \textbf{84}, 119--135
(2000).


\bibitem{Arnold1} V. I. Arnold, Bernoulli--Euler up-down numbers associated
with function singularities, their combinatorics and arithmetics, Duke Math.
J. \textbf{63}, 537--554 (1991).

\bibitem{Arnold2} V. I. Arnold, The calculus of snakes and the combinatorics
of Bernoulli, Euler and Springer numbers of Coxeter groups, Russian Math.
Surveys, \textbf{47}, 1--51 (1992).

\bibitem{AB} A. Baker, Combinatorial and arithmetic identities based on
formal group laws, Lect. Notes in Math. \textbf{1298}, 17--34, Springer
(1987).

\bibitem{BCRS} A. Baker, F. Clarke, N. Ray, L. Schwartz, On the Kummer
congruences and the stable homotopy of BU, Trans. Amer. Math. Soc. \textbf{%
316}, 385--432 (1989).


\bibitem{Bernoulli} J. Bernoulli, \textit{Ars conjectandi}, Basel, 1713.

\bibitem{BB} R. P. Boas, R. C. Buck, \textit{Polynomial Expansions of
Analytic Functions}, Springer--Verlag, New York, 1964.

\bibitem{BL} A. di Bucchianico and D. Loeb, Umbral Calculus, The Electronic
Journal of Combinatorics, \textbf{DS3} (2000), http://www.combinatorics.org

\bibitem{Carl} L. Carlitz, Degenerate Stirling, Bernoulli and Eulerian
numbers, Utilitas Math. \textbf{15}, 51--88 (1979).

\bibitem{Carlitz} L. Carlitz, Some Theorems on Bernoulli numbers of higher
order, Pacific J. Math. \textbf{2}, 127--139 (1952).

\bibitem{CO} L. Carlitz and F. R. Olson, Some Theorems on Bernoulli and
Euler numbers of higher order, Duke Math. J. \textbf{21}, 405--421 (1954).


\bibitem{Clarke} F. Clarke, The universal von Staudt Theorems, Trans. Amer.
Math. Soc. \textbf{315}, 591--603 (1989).


\bibitem{DLM} B. Doyon, J. Lepowsky, A. Milas, Twisted vertex operators and
Bernoulli polynomials, arXiv: math. QA/0311151 (2003).

\bibitem{DSS} K. Dilcher, L. Skula and I. SH. Slavutskii, \textit{Bernoulli
numbers. Bibliography (1713--1990)}, Queen's Papers in Pure and Appl. Math.
\textbf{87}, Queen's University, Kingston, ON, 1991, iv+175 pp. (1990),
http://www.mathstat.dal.ca/dilcher/$\sim $bernoulli.html

\bibitem{DS} G. V. Dunne, C. Schubert, Bernoulli Number Identities from
Quantum Field Theory, arXiv: math. NT/0406610 (2004).

\bibitem{Eie} M. Eie, A note on Bernoulli numbers and Shintani generalized
Bernoulli polynomials, Trans. Amer. Math. Soc. \textbf{348}, 1117--1136
(1996).

\bibitem{Erdelyi} A. Erd\'{e}lyi, \textit{The Bateman Manuscript Project},
vol 2, McGraw--Hill, 1953.

\bibitem{FV} D. B. Fairlie and A. P. Veselov, Faulhaber and Bernoulli
polynomials and solitons, Physica D \textbf{152--153}, 47--50 (2001).

\bibitem{FP} C. Faber and R. Pandharipande, Hodge integrals and
Gromov--Witten theory, Invent. Math. \textbf{139}, 173--199 (2000).

\bibitem{GV} M. P. Grosset and P. Veselov, Lam\'{e} equation, quantum top
and elliptic Bernoulli polynomials, arXiv: math--ph/0508068 (2005).

\bibitem{Haze} M. Hazewinkel, \textit{Formal Groups and Applications},
Academic Press, New York, 1978

\bibitem{How} F. T. Howard, Congruences and recurrences for Bernoulli
numbers of higher order, Fibonacci Quart. \textbf{32}, 316--328 (1994).

\bibitem{IR} K. Ireland and M. Rosen, \textit{A Classical Introduction to
Modern Number Theory}. Second edition. Springer--Verlag, 1990.

\bibitem{Iwa} K. Iwasawa, \textit{Lectures on p--adic L--Functions},
Princeton University Press, Princeton, 1972

\bibitem{Koblitz} N. Koblitz, \textit{p--adic Numbers, p--adic Analysis and
Zeta--functions}, Springer--Verlag, New York (2nd edition), 1984.

\bibitem{Lang2} S. Lang, \textit{Introduction to Modular Forms},
Springer--Verlag, 1976.

\bibitem{LNO} D. Levi, J. Negro and M. A. del Olmo, Discrete derivatives and
symmetries of difference equations, J. Phys. A: Math. Gen. \textbf{34},
2023--2030 (2001).

\bibitem{LTW1} D. Levi, P. Tempesta and P. Winternitz, Umbral Calculus,
Diffference Equations and the Discrete Schr\"{o}dinger Equation\textit{, }J.
Math. Phys. \textbf{45}, 4077--4105 (2004).

\bibitem{LTW2} D. Levi, P. Tempesta and P. Winternitz, Lorentz and Galilei
Invariance on Lattices, Phys. Rev. D \textbf{69}, 105011 (2004).

\bibitem{MT} S. Marmi, P. Tempesta, Polylogarithms, hyperfunctions and generalized Lipschitz summation formulae, Preprint I-2007, Centro di ricerca matematica Ennio De Giorgi, Scuola Normale Superiore, Pisa, 2007. 

\bibitem{Mil} H. Miller, Universal Bernoulli numbers and the S$^{1}$ --
transfer, \textit{Current Trends in Algebraic Topology 2}, pt. 2, 437--449,
CMS--AMS (1982).

\bibitem{Phil} G. M. Phillips, \textit{Interpolation and Approximation by
Polynomials}, CMS Books in Mathematics, Springer--Verlag, 2003.

\bibitem{Ray} N. Ray, Stirling and Bernoulli numbers for complex oriented
homology theory, in Algebraic Topology, Lecture Notes in Math. 1370, pp.
362--373, G. Carlsson, R. L. Cohen, H. R. Miller and D. C. Ravenel (Eds.),
Springer--Verlag, 1986.

\bibitem{Riordan1} J. Riordan, \textit{An Introduction to Combinatorial
Analysis}, Wiley, New York, 1958.

\bibitem{Roman} S. Roman, \textit{The Umbral Calculus}, Academic Press, New
York, 1984.

\bibitem{Rota} G. C. Rota, \textit{Finite Operator Calculus}, Academic
Press, New York, 1975.

\bibitem{Slav1} I. Sh. Slavutskii, Outline of the history of research on the
arithmetic properties of Bernoulli numbers, Istor.--Mat. Issled. \textbf{%
32/33}, 158-- 181 (1990) (Russian).

\bibitem{T1} Y. Smirnov and A. Turbiner, Lie algebraic discretization of
differential operators, Mod. Phys. Lett A \textbf{10}, 1795--1802 (1995).

\bibitem{Tempesta1} P. Tempesta, Formal Groups, Bernoulli--type polynomials
and L--series, C. R. Math. Acad. Sci. Paris, to appear.

\bibitem{Tempesta2} P. Tempesta, Dirichlet and Riemann Hurwitz zeta
functions associated with formal groups and generalized Bernoulli
polynomials, preprint (2007)

\bibitem{T2} A. Turbiner, Canonical Discretization I. Different faces of the
(an)harmonic oscillator, Int. J. Mod. Phys. A \textbf{16}, 1579--1603 (2001).


\bibitem{Wash} L. C. Washington, \textit{Introduction to Cyclotomic Fields}.

\bibitem{You} P. T. Young, Congruences for Bernoulli, Euler and Stirling
Numbers, J. Number Theory \textbf{78}, 204--227 (1999).

Second edition. Springer--Verlag, 1997.
\bibitem{Zag} D. Zagier, Valeurs des foctions z\^eta des corps quadratiques r\'eels aux entiers n\'egatifs, Ast\'erisque 41-42, 135-151 (1977).
\end{thebibliography}
\end{document}